\documentclass{amsart}
\usepackage{amsthm,amsfonts,amsmath,amssymb,latexsym,epsfig,upref,eucal,ae}
\usepackage{mathrsfs}  
\usepackage{tikz-cd}

\usepackage{cite}
\usepackage{scalerel,stackengine}
\usepackage{textcomp}

\usepackage{mathtools}

\usepackage{hyperref}

\theoremstyle{plain}
\newtheorem{theorem}{Theorem}[section]

\theoremstyle{definition}

\newtheorem{definition-theorem}[theorem]{Definition-Theorem}

\theoremstyle{remark}
\newtheorem{remark}[theorem]{Remark}

\def\cF{\mathcal{F}}

\def\cO{\mathcal{O}}

\def\cT{\mathcal{T}}

\def\R{\mathbb{R}}

\def\Z{\mathbb{Z}}
\def\N{\mathbb{N}}

\def\C{\mathbb{C}}

\def\Om{\Omega}

\def\>{\rangle}

\def\<{\langle}
\def\>{\rangle}

\def\Hom{\mathrm{Hom}}
\def\Spec{\mathrm{Spec}}

\begin{document}

\title[Zariski--Lipman conjecture for toric varieties]
{The Zariski--Lipman conjecture for toric varieties}

\author[C. Tipler]{Carl Tipler}
\address{Univ Brest, UMR CNRS 6205, Laboratoire de Math\'ematiques de Bretagne Atlantique, France}
\email{carl.tipler@univ-brest.fr}

\date{\today}

 \begin{abstract}
We give a short proof of the Zariski--Lipman conjecture for toric varieties : any complex toric variety with locally free tangent sheaf is smooth.
\end{abstract}

\maketitle

 
\section{Introduction}
\label{sec:intro}
The Zariski--Lipman conjecture states that an algebraic variety $X$ over $\C$ with locally free tangent sheaf is necessarily smooth. This conjecture has been proved under various additional assumptions, see for example \cite{Hoch} in the $\N$-graded case, \cite{Kall} for complete intersections, \cite{Druel} for log canonical varieties and \cite{Graf,BerGra,Bis} for surfaces. In this note, we give a direct and simple proof for toric varieties, that is for irreducible varieties endowed with an algebraic effective action of a complex torus with a dense open orbit. Recall that the tangent sheaf $\cT_X$ of an algebraic variety $X$ is the dual of its sheaf of K\"ahler differentials, that is $\cT_X=\mathrm{Hom}_{\cO_X}(\Omega^1_X,\cO_X)$.
\begin{theorem}
\label{theo:intro}
 Let $X$ be a toric variety over $\C$. Assume that the tangent sheaf of $X$ is locally free. Then $X$ is smooth.
\end{theorem}
As we will see in Remark \ref{rem:normal} and Remark \ref{rem:nakai}, Theorem \ref{theo:intro} can be obtained as a corollary of the results from \cite{GrafKov} or from \cite{LevStaf}. Nevertheless, the proof presented here is fairly simple and relies purely on methods from toric geometry.
\begin{remark}
\label{rem:normal}
 If $(X,0)$ is a germ of a rational singularity over $\C$, then $(X,0)$ is Cohen--Macaulay. If moreover $(X,0)$ has a locally free tangent sheaf, its canonical divisor is Cartier, and thus $(X,0)$ is Gorenstein. As noticed in \cite{Graf}, rational Gorenstein singularities are canonical, a class of singularities for which the Zariski--Lipman conjecture has been proved in \cite{GrafKov,Druel}. Normal toric singularities being rational \cite[Theorem 11.4.2]{CLS}, Theorem \ref{theo:intro} is a corollary of the results in \cite{GrafKov}. 
\end{remark}

\begin{remark}\label{rem:nakai}
Locally, a normal toric variety with no torus factor is a categorical quotient of a smooth affine toric variety by a reductive group action \cite[Theorem 5.1.11]{CLS}. From \cite[Corollary 5.11]{LevStaf}, this implies that the Nakai conjecture holds true for normal toric varieties. As the Nakai conjecture implies the Zariski-Lipman conjecture \cite[Proposition 2]{Rego}, Theorem \ref{theo:intro} follows from \cite{LevStaf}. The author would like to thank Thierry Levasseur for pointing to him this reference.
\end{remark}

\begin{remark}
The Zariski--Lipman conjecture is actually stated over any field of characteristic zero. As noticed in \cite{Hoch}, with no loss of generality one may assume the field to be algebraically closed. In this note, we will work over the complex numbers, but Theorem \ref{theo:intro} and its proof hold over any algebraically closed field of characteristic zero.
\end{remark}

\subsection*{Acknowledgments}  
The author would like to thank Patrick Graf for explaining Remark \ref{rem:normal} and Thierry Levasseur for explaining Remark \ref{rem:nakai}, as well as Thibaut Delcroix, Henri Guenencia and the anonymous referee for useful remarks. The author is partially supported by the grants MARGE ANR-21-CE40-0011 and BRIDGES ANR--FAPESP ANR-21-CE40-0017.


\section{Klyachko's description of toric reflexive sheaves}
\label{sec:toricbasics}
Let $X$ be a toric variety of dimension $n$ over $\C$, with torus $T_N=N\otimes_\Z \C^*$, for $N$ the rank $n$ lattice of its one-parameter subgroups. Denote by $M=\Hom_\Z(N,\Z)$ its character lattice. Let $\cT_X=\mathrm{Hom}_{\cO_X}(\Omega^1_X,\cO_X)$  be the tangent sheaf of $X$, that is the dual of its sheaf of K\"ahler differentials (see e.g. \cite[Section 8.0]{CLS}). According to \cite[Theorem 3]{Lip}, if $\cT_X$ is locally free, then $X$ is normal, which we will assume from now on. Then, $X$ is the toric variety associated to a fan $\Sigma$ of strongly convex rational polyhedral cones in $N_\R=N\otimes_\Z \R$ \cite[Chapter 3]{CLS}. In particular, $X$ is covered by the $T_N$-invariant affine varieties $U_\sigma=\Spec(\C[M\cap\sigma^\vee])$, for $\sigma\in \Sigma$.

Recall that a coherent sheaf $\cF$ on $X$ is called $T_N$-equivariant if there is an isomorphism
$\varphi : \alpha^*\cF \to \pi_2^*\cF$
satisfying some cocycle condition (see e.g. \cite[Section 5]{perl} or \cite{Her}) where $\alpha : T_N \times X \to X$,  $\pi_1 : T_N\times X \to T_N$ and $\pi_2 : T_N \times X \to X$ stand for respectively the $T_N$-action, the projection on $T_N$ and the projection on $X$. The $T_N$-action on $X$ induces naturally an equivariant structure on its sheaf of K\"ahler differentials given by the composition :
$$
\alpha^*\Om^1_X \xrightarrow{d\alpha} \Om^1_{T_N\times X}\xrightarrow{\simeq} \pi_1^*\Om^1_{T_N}\oplus\pi_2^*\Om^1_X \xrightarrow{\mathrm{pr}_2} \pi_2^*\Om^1_X
$$
where $\mathrm{pr}_2$ is the projection on the second factor \cite[Section 2]{Her}.
By duality, $\cT_X$ is $T_N$-equivariant, and, being dual to the coherent sheaf $\Om^1_X$, it is reflexive. Klyachko proved that equivariant reflexive sheaves on toric varieties are described by families of filtrations \cite{Kl} (see also \cite{perl}). Let us recall briefly this description for the tangent sheaf. As $\cT_X$ is reflexive, its sections extend over codimension $2$ subvarieties. Denote by $\Sigma(l)$ the set of $l$-dimensional cones in $\Sigma$, and consider
$$
X_0=\bigcup_{\sigma\in\Sigma(0)\cup\Sigma(1)} U_\sigma.
$$ 
By the orbit-cone correspondence \cite[Section 3.2]{CLS}, $X_0$ 
is the complement of $T_N$-orbits of co-dimension greater or equal to $2$, and thus $\cT_X=\iota_*\cT_{X_0}$, where $\iota : X_0 \to X$ is the inclusion. By equivariance, $\cT_{X_0}$ is entirely characterised by the sections $\Gamma(U_\sigma,\cT_X)$, for $\sigma\in\Sigma(0)\cup\Sigma(1)$. If $\sigma=\lbrace 0 \rbrace$, $U_{\lbrace 0\rbrace}=T_N$, and
$$
\Gamma(U_{\lbrace 0\rbrace},\cT_X)=N_\C\otimes_\C \C[M],
$$
for $N_\C:=N\otimes_\Z\C$ is the Lie algebra of $T_N$. Then, if $\rho\in \Sigma$ is a ray (i.e. a one-dimensional cone), $\Gamma(U_{\rho}, \cT_X)$ is graded by $M/( M\cap \rho^\perp)\simeq\Z$. As $T_N$ is a dense open subset of $U_\rho$, the restriction map $\Gamma(U_\rho,\cT_X)\to\Gamma(U_{\lbrace 0\rbrace},\cT_X)$ is injective and induces a decreasing $\Z$-filtration 
$$
  \ldots \subset E^\rho(i)\subset E^\rho(i-1) \subset \ldots \subset N_\C
$$
such that one has
 \begin{equation*}
  \label{eq:sheaf from family of filtrations}
  \Gamma(U_{\rho}, \cT_X)=\bigoplus_{m\in M} E^\rho(-\langle m,u_\rho\rangle)\otimes \chi^m,
 \end{equation*}
 where we denote by $u_\rho$ the primitive generator of $\rho$ and $\langle \cdot,\cdot\rangle$ the duality pairing.
Explicitly, the family of filtrations $(E^\rho(\bullet))_{\rho\in\Sigma(1)}$ for $\cT_X$ is given by \cite[Example 2.3(5) on page 350]{Kl}:
$$
E^\rho(i)=\left\{ \begin{array}{ccc}
                  N_\C & \textrm{ if } & i\leq 0 \\
                  \C\cdot u_\rho & \textrm{ if } & i= 1 \\
                  \lbrace 0 \rbrace & \textrm{ if } & i\geq 2.
                 \end{array} 
\right.
$$
Finally, by Klyachko's compatibility condition \cite[Theorem 2.2.1]{Kl} (see \cite[Section 5]{perl} for a detailed treatment on normal toric varieties), $\cT_X$ is locally free if and only if the family of filtrations $(E^\rho(\bullet))_{\rho\in\Sigma(1)}$ satisfies that for each $\sigma\in\Sigma$, there exists a decomposition 
$$
N_\C= \bigoplus_{[m]\in M/(M\cap\sigma^\perp)} E^\sigma_{[m]}
$$
such that for each ray  
$\rho\subset \sigma$:
$$
E^\rho(i)=\bigoplus_{\langle m,u_\rho\rangle\geq i}  E^\sigma_{[m]}.
$$
\section{Proof of Theorem \ref{theo:intro}}
\label{sec:proof}
Assume from now on that $\cT_X$ is locally free. We want to show that $X$ is smooth. As this is a local condition, we might as well assume $X$ affine, so that $X=U_\sigma$ for some strongly convex rational polyhedral cone $\sigma\subset N_\R$. With no loss of generality, we can also assume that $X$ has no torus factor, so that $\lbrace u_\rho,\, \rho\in\sigma(1)\rbrace$ spans $N_\R$ \cite[Proposition 3.3.9]{CLS}. Then, $X$ is smooth if and only if $\sigma$ is smooth \cite[Theorem 1.3.12]{CLS}, which, by definition, is equivalent to the fact that $\lbrace u_\rho, \rho\in\sigma(1)\rbrace$ is a $\Z$-basis for $N$.

As $\lbrace u_\rho,\, \rho\in\sigma(1)\rbrace$ spans $N_\R$,  
 $\sigma^\perp=\lbrace 0 \rbrace$. Then by Klyachko's compatibility condition, we can find a decomposition 
$$
N_\C= \bigoplus_{m\in M} E^\sigma_{m}
$$
with, for $\rho\in\sigma(1)$, and $i\in\lbrace 2, 1, 0 \rbrace$,
\begin{equation}
 \label{eq:Klyachkos}
E^\rho(i)=\bigoplus_{\langle m,u_\rho\rangle\geq i}  E^\sigma_{m}.
\end{equation}
First, from $E^\rho(2)=\lbrace 0 \rbrace$, we deduce that $E^\sigma_m\neq\lbrace 0 \rbrace$ only if for all $\rho$, $\langle m,u_\rho\rangle\leq 1 $. 
Secondly, taking $i=1$ in (\ref{eq:Klyachkos}),
$$
\C\cdot u_\rho=\bigoplus_{\langle m,u_\rho\rangle\geq 1}  E^\sigma_{m}.
$$
Thus, for each $\rho\in \sigma(1)$ we can find $m_\rho\in M$ with $\langle m_\rho,u_\rho\rangle=1$ and such that $\C \cdot u_\rho = E_{m_\rho}^\sigma$. Note that $\sigma$ contains no line by strong convexity, so if $\rho \neq \rho'$, then $u_\rho \notin \Z\cdot u_{\rho'}$ and $\C\cdot u_\rho\neq \C\cdot u_{\rho'}$. As $\lbrace u_\rho, \rho\in\sigma(1)\rbrace$ spans $N_\C$ over $\C$, 
$$
N_\C  =  \sum_{\rho\in \sigma(1)}\C \cdot u_\rho
     =  \bigoplus_{\rho\in \sigma(1)} E^\sigma_{m_\rho} 
       \subset  \bigoplus_{m\in M} E^\sigma_{m}=N_\C.
$$
Then, 
$$N_\C=  \bigoplus_{\rho\in \sigma(1)} E^\sigma_{m_\rho}$$
and $\sigma(1)$ contains exactly $n$ elements. Let $\rho,\rho'\in\sigma(1)$ be two distinct rays. Necessarily, 
$$\C\cdot u_\rho\cap \C\cdot u_{\rho'}=\lbrace 0\rbrace$$
and by (\ref{eq:Klyachkos}) with $i=1$, $m_\rho$ must satisfy $\langle m_\rho, u_{\rho'}\rangle \leq 0$. 

Last, taking now $i=0$ in (\ref{eq:Klyachkos}), we deduce from $E^\rho(0)=N_\C$ that $\langle m_\rho, u_{\rho'} \rangle = 0$. 
To conclude, for all $\rho, \rho'\in\sigma(1)$,
$$
\langle m_\rho, u_{\rho'}\rangle =\left\{
\begin{array}{ccc} 
1 & \mathrm{ if } & \rho=\rho'\\
0 & \mathrm{ if } & \rho\neq\rho'.
\end{array}
\right.
$$ 
Hence, each element $u\in N$ can be uniquely written 
$$
u=\sum_{\rho\in\sigma(1)} \langle m_\rho , u \rangle\, u_\rho
$$
with $\langle m_\rho , u \rangle\in\Z$ for each $\rho\in\sigma(1)$. Thus, $\lbrace u_\rho, \rho\in\sigma(1)\rbrace$ is a basis of $N$, which ends the proof of Theorem \ref{theo:intro}.

\end{document}